\documentclass[letterpaper, 10 pt, conference]{ieeeconf}  % Comment this line out
 
\IEEEoverridecommandlockouts                              % This command is only

\usepackage{graphics} % for pdf, bitmapped graphics files

\usepackage{epsfig} % for postscript graphics files

\usepackage{mathptmx} % assumes new font selection scheme installed

\usepackage{times} % assumes new font selection scheme installed
 
\usepackage{amsmath} % assumes amsmath package installed

\usepackage{amssymb}  % assumes amsmath package installed

\usepackage{color}

\usepackage{algorithm}

\usepackage{algorithmic}

\newtheorem{thm}{Theorem}

\newtheorem{dfn}{Definition}

\newtheorem{prop}{Proposition}

\newcommand{\RR}{{\mathbb R}}

\newcommand{\norm}[1]{\left| \left| #1 \right| \right|}

\newcommand{\dotex}{{\frac{d}{dt}}}

\title{\LARGE \bf
An intrinsic Cram{\'e}r-Rao bound on SO(3) for (dynamic) attitude filtering}

\author{Silv\`ere Bonnabel and Axel Barrau% <-this % stops a space
\thanks{A. Barrau and S. Bonnabel are with MINES ParisTech, PSL Research University, Centre for robotics, 60 Bd St Michel 75006 Paris, France
       {\tt\footnotesize [axel.barrau,silvere.bonnabel]@mines-paristech.fr}}%
}

\include{myLatex}

\graphicspath{{Figures/}}

\begin{document}

\maketitle

\thispagestyle{empty}

\pagestyle{empty}

%%%%%%%%%%%%%%%%%%%%%%%%%%%%%%%%%%%%%%%%%%%%%%%%%%%%%%%%%%%%%%%%%%%%%%%%%%%%%%%%

\begin{abstract}
In this note an intrinsic version of the Cram{\'e}r-Rao bound on estimation accuracy is  established on the Special Orthogonal group $SO(3)$. It is intrinsic in the sense that it does not rely on a specific choice of coordinates on $SO(3)$: the result is derived using rotation matrices, but remains valid when using other parameterizations, such as quaternions. For any estimator $\hat R$ of $R\in SO(3)$ we give indeed a lower bound on the covariance matrix of $ \log \left( R\hat R^T \right)$, that is, the estimation error expressed in terms of group multiplication, whereas the usual estimation error $ \left( \hat R-R \right)$ is meaningless on $SO(3)$.  The result is first applied to Whaba's problem. Then, we consider the problem of a continuous-time nonlinear \emph{deterministic} system on $SO(3)$ with discrete measurements subject to additive isotropic \emph{Gaussian noise}, and we derive a lower bound to the estimation error covariance matrix. We prove the \emph{intrinsic} Cram{\'e}r-Rao bound coincides with the covariance matrix returned by the \emph{Invariant} EKF, and thus \emph{can} be computed online. This is in sharp contrast with the general case, where the bound can only be computed if the true trajectory of the system is known. 
\end{abstract}
\section{Introduction}

We consider the problem of attitude estimation from vector measurements, where the attitude parameter is static in the first place, and then the dynamic case where it evolves on the special orthogonal group $SO(3)$. Attitude estimation, both in the static and dynamic cases, has been the subject of numerous works due to its potential applications to e.g.  aerial vehicles (or satellites) control. The references are too numerous to be exhaustively listed, and the reader is referred to e.g.   \cite{madinehi2013rigid} for an overview of estimation problems on SO(3), or the survey  \cite{crassidis2007survey}, and e.g. the paper  \cite{izadi2014rigid} for a very recent work on the subject. 

In the present paper an intrinsic version of the Cram{\'e}r-Rao bound on estimation accuracy is   established in Section \ref{smith:sec} on the Special Orthogonal group $SO(3)$. It is intrinsic in the sense that it does not rely on a specific choice of coordinates on $SO(3)$. For any estimator $\hat R$ of $R\in SO(3)$ we give a lower bound on the covariance matrix of $ \log \left( R\hat R^T \right) $, that is, the estimation error expressed in terms of group multiplication, and then projected onto a three dimensional vector space using the logarithmic map of $SO(3)$. This error indeed makes sense as $R\hat R^T$  is the rotation that maps the estimated orientation to the true orientation, whereas the usual estimation error $\left( \hat R-R \right)$ is meaningless on $SO(3)$, as $\left( \hat R-R \right)$ is not a rotation matrix, and has no intrinsic counterpart. 
Taking advantage of the Lie group structure of the space the calculations are rather simple and direct. Viewing $SO(3)$ as a manifold and choosing an invariant metric, we  recover in a simple way the result derived by S. Smith in \cite{smith-2005} (see also the recent work of N. Boumal \cite{boumal2013intrinsic,boumal2014cramer}).

As a straightforward application, the result is first applied to the static attitude estimation problem, also known as Whaba's problem \cite{wahba1965}, for which we derive a lower bound. Note a (classical) Cram{\'e}r-Rao lower bound has already been proposed  in \cite{crassidis2000optimal} for the linearized problem.  Then, we consider in Section \ref{taylor:sec} the  problem of attitude estimation/filtering from vector measurements and angular velocity measurements from a gyroscope (see e.g. \cite{izadi2014rigid,mahony2009nonlinear} and \cite{hua2014implementation} for an implementation)  in the degenerate case where the gyroscope is of much higher quality than the other sensors. For systems possessing deterministic dynamics and stochastic output measurements, J. H. Taylor proved in \cite{taylor1979} the Cram{\'e}r-Rao bound is provided by the Extended Kalman Filter (EKF) covariance, linearized around the \emph{true} unknown trajectory of the system, and thus cannot be computed online. Thanks to the invariance properties of the system, we prove  the Cram{\'e}r-Rao bound does not depend on the true system's trajectory, and can be computed online.

The Invariant Extended Kalman Filter (IEKF) is a recent methodology to modify the EKF in order to account for the invariance properties of the state space when devising EKFs on Lie groups, see \cite{bonnabel2007left,bonnabel2009invariant}, and more recently \cite{barrau2013intrinsicp} where an IEKF is derived on $SO(3)$ with discrete time observations. A remarkable property of the IEKF, akin to the properties of symmetry-preserving observers \cite{bonnabel2008symmetry,arxiv-08} from which the IEKF is derived,   is that the estimation error system depends on the system's trajectory in a reduced manner, and sometimes does not depend on it at all, a property shared by the intrinsic Cram{\'e}r-Rao bound derived in this paper. In fact, the links between both theories go beyond: in the case considered here on $SO(3)$, we prove in Section \ref{kalman:sec} the \emph{intrinsic} Cram{\'e}r-Rao bound coincides with the covariance matrix returned by the \emph{Invariant} EKF (and thus can be computed online).

\section{An intrinsic Cram{\'e}r-Rao bound on SO(3)}\label{smith:sec}

\subsection{The classical Cram\'er-Rao bound}

Consider a family of probability densities   $p \left( x \mid \theta \right)$  parameterized by a vector $\theta\in\RR^k$. Consider an unbiased estimator $\hat\theta(x)$ of the parameter $\theta$ from a sample measurement $x$. The requirement that $\hat\theta$ be unbiased means it must be ``good'' (for a large sample) whatever $\theta$. Because of this requirement, the estimator can not recover $\theta$ exactly, given a single or a finite number of measurements. This fact is formalized by the well-known existence of a lower bound on the accuracy of the estimator: the so-called Cram\'er-Rao bound. Mathematically, it states the average estimation error covariance is lower bounded as follows
$$
P := \textbf{E}_{\theta} \left(  \left( \hat\theta-\theta \right) \left( \hat\theta-\theta \right)^T \right) \succeq  J^{-1}
$$where $ \textbf{E}_{\theta}$ denotes the expectation with respect to the probability law  $p \left( x \mid \theta \right)$, and the matrix $J$ (the inequality is in the sense of the Loewner order) is the so-called Fisher information matrix, defined as the Hessian with respect to $\theta$ of the average log likelihood $\textbf{E}_{\theta} \left( \ln \left( p \left(x\mid\theta \right) \right) \right)$, where $\ln(.)$ denotes the natural logatithm of $\mathbb{R}^+$. 

The interesting question raised by S.T. Smith \cite{smith-2005}, is whether there exists an analogue of this bound for a parameter $\theta$ that belongs to a Riemannian manifold, and not a vector space anymore. This kind of question can arise in signal processing, where one seeks to estimate for instance a subspace, as in Principal Component Analysis, that is, an element of the Grassman manifold. To answer this question, one must first find a way to compare the estimator $\hat \theta$ and the true parameter $\theta$, as on a manifold the quantity $\hat\theta-\theta$ has no meaning. This can be done through the Riemann exponential map, and then, \cite{smith-2005} proves a Cram\'er-Rao bound can be produced. Adapting the classical proof to the manifold case, he shows that the (well-defined) error covariance $P$ is lower bounded by a (well-defined) information matrix $J$, plus additional terms steming from the curvature of the parameter space. Unfortunately, the formula is not in closed form. However, for sufficiently small covariance $P$ it can be expanded up to terms of order $O(P^{3/2})$. 

The quite inspiring paper \cite{smith-2005} draws new links between statistics and geometry. It has been in particular adapted to the quotient manifold case in \cite{boumal2013intrinsic}. In the present paper, we derive similar results for the $SO(3)$ case, more simply, and we apply them to two attitude estimation problems. Note that, a tutorial presentation on the intrinsic Cram{\'e}r-Rao bound can be found in e.g. \cite{barrau2013note}.

\subsection{Direct derivation of the Cram\'er-Rao bound  on $SO(3)$ }

We compute here the Intrinsic Cram{\'e}r-Rao Lower Bound (ICRLB) on $SO(3)$, up to the second order terms in the estimation error $d(\hat R,R)$ where $d$ denotes the bi-invariant distance on $SO(3)$. This allows to recover in a simple and direct way the result of \cite{smith-2005} taking advantage of the Lie group structure of the parameter space (and thus without having to evaluate Riemann's curvature tensor at $R$).

\subsubsection{Preliminaries}

$SO(3)$ is a Lie group of dimension $3$, and thus a Riemannian manifold. The tangent space  at $Id$, the Identity rotation, denoted $\mathfrak{so(3)}$ is called the Lie Algebra of $SO(3)$ and can be identified with $\RR^3$ that is
$$
\mathfrak{so(3)}\approx \RR^3
$$
 Using rotation matrices (i.e. viewing $SO(3)$ as a submanifold of $\RR^{3\times 3}$)  the (group) exponential map defined by 
\begin{align*}
&\exp : \mathfrak{so(3)}\mapsto SO(3)
\\ &\exp(\xi)=\text{expm} [(\xi)_\times]
\end{align*}
where expm denotes the matrix exponential map, and where $(a)_\times\in\RR^{3\times 3}$ for $a\in\RR^3$ denotes the skew symmetric matrix defined by $(a)_\times u=a\times u,~\forall u\in\RR^3$. In a neighborhood of $Id$, the exponential map can be inverted. The (group) logarithmic  map
$$
\log : SO(3)\mapsto \mathfrak{so(3)}
$$
is defined as the inverse of exp. For any estimator $\hat R$ of a parameter $R\in SO(3)$, it allows to measure the mean quadratic estimation error
$$
\textbf{E}_R \left( \norm{\log \left(R\hat R^T \right) }_{\RR^3}^2 \right)= \textbf{E}_R\left(d(R,\hat R)^2 \right)
$$
 The logarithmic map allows also to define a covariance matrix of  the (right-invariant) estimation error 
\begin{align}
P = \textbf{E}_R \left( \log \left(R\hat R^T \right) \log \left( R\hat R^T \right)^T \right) \in\RR^{3\times 3}\label{P:def}
\end{align}
and we have as usual $\textbf{E}_R \left( d \left( R,\hat R \right)^2 \right) = Tr \left( P \right)$. 

Note that, if $R$ and $\hat R$ denote orientations of bodies in space, the quantity $\log \left( R\hat R^T \right) \in\RR^3$  has an intrinsic purely geometrical interpretation as its orientation denotes the axis of rotation and its norm the rotation time around which the body with configuration $\hat R$ must rotate in order to reach configuration $R$, and provides us with a natural orientation error in $\RR^3$.

\subsubsection{Main result} 
Consider the family of densities parameterized by an element $R$  of SO(3)
\[
p \left( X \mid R \right), X\in \mathcal M,~R \in SO(3)
\]where the sample space $\mathcal M$ is a measurable space. To fix ideas we will consider in the sequel that  $\mathcal M$ is $\RR^q$ for $q \in \mathbb{N}$, as will be the case in the examples. Using the exponential map, and mimicking its Riemannian analogue \cite{smith-2005}, we introduce
\begin{dfn} The intrinsic Fisher information matrix can be defined in a right-invariant basis  as follows for any $\xi\in\RR^3$
\begin{align}
\xi^TJ(R)\xi & = \int \left( \left. \dotex \right|_{t=0} \ln p \left(X \mid \exp(t\xi)R \right) \right) \ldots  \nonumber\\
& \left( \left. \dotex \right|_{t=0} \ln p \left( X \mid \exp(t\xi)R \right) \right)^Tp(X|R)dX
\end{align}and then $J(R)$ can be recovered using the standard polarization formulas:
$$
\xi^TJ\nu=\frac{1}{2} \left( \left(\xi+\nu \right)^TJ \left( \xi+\nu \right)-\xi^TJ\xi-\nu^TJ\nu \right).
$$
\end{dfn}
Besides, note that, using the fact that 
\begin{align}
\int p \left( X\mid \exp \left( t\xi \right)R \right) & \dotex\ln  p \left( X\mid \exp \left( t\xi \right) R \right) dX \nonumber \\
& = \dotex \int p(X\mid \exp \left( t\xi \right)R)dX = 0
\end{align}
 and differentiating the latter equality w.r.t $t$ then reusing that $\dotex p=p\dotex \log p$ we have
\begin{align*}
\int & \left( \dotex\ln p \left( X\mid \exp \left( t\xi \right) R \right) \right) \left(\dotex\ln p \left( X\mid \exp \left( t\xi \right) R \right) \right) dX \\
 & + \int p \left( X\mid \exp \left(t\xi \right) R \right) \frac{d^2}{dt^2}\ln p \left( X\mid \exp \left(t\xi\right)R \right) dX = 0
\end{align*}
allowing to recover an intrinsic version of the classical result according to which the information matrix can be also defined using a second order derivative (i.e. a Hessian). 
\begin{prop}The intrinsic Fisher information matrix also writes$$
{\xi^TJ(R)\xi=- \textbf{E}_R \left( \left. \frac{d^2}{dt^2} \right|_{t=0} \ln p \left( X\mid \exp \left( t\xi \right)R \right) \right)}
$$\end{prop}
Let $\hat{R}$ be an unbiased estimator of $R$ in the sense of centered intrinsic (right invariant) error $R \hat R^T$, that is,
\[
\int_X \log \left[ R \hat R^T \left( X \right) \right] p \left(X \mid R \right) dX =  0.
\]
\begin{thm}Let $P$ be the covariance matrix of the estimation error as defined in \eqref{P:def}. We have then
\begin{equation}
\begin{aligned}
\label{result2}
P \succeq  J \left( R \right)^{-1}  -&\frac{1}{12} \left( Tr \left( P \right) I_3 - P \right) J \left( R \right)^{-1} \\&\quad -\frac{1}{12} J \left( R \right)^{-1} \left( Tr \left( P \right) I_3 - P \right) \end{aligned}
\end{equation}
where we have neglected terms of order $\textbf{E}_R \left(||\log[R \hat R^T (X) ]||^3\right)$. 
For small errors, we can neglect the terms in $P$ on the right hand side (curvature terms) yielding
\begin{align*}
P & =\int_X \log \left( R\hat R \left( X \right)^T \right) \log \left( R\hat R \left( X \right)^T \right)^T p \left( X \mid R \right) dX \\
 & \succeq J(R)^{-1} \text{+ C},
\end{align*}
where C are terms of higher order, linked to the effect of the curvature of the parameter space $SO(3)$, hence the letter C, which here stands for ``curvature terms".\end{thm}

\subsection{Application to Wahba's problem}
Before proving the results let us give an example application. In this subsection we assume measurements are of the form
$$
X_k=R^Td_k+V_k,
$$
where $d_k$'s are some reference vector in $\RR^3$, and where $V_k$'s are independent  isotropic and Gaussian noises with covariance matrices $\sigma_k^2 Id$. This implies the following density form
%\begin{align*}
%- \log p(X_k|R) & = \text{ Cste }+\frac{1}{2\sigma_k^2}(X-R^Td_k)^T(X-R^Td_k) \\
% & = \frac{1}{2\sigma_k^2}\norm{X-R^Td_k}^2
%\end{align*}
\begin{align*}
- \ln p \left( X_k \mid R \right) = \frac{1}{2\sigma_k^2}\norm{X-R^Td_k}^2.
\end{align*}
Wahba's problem consists in finding the maximum likelihood estimator of $R$ and has been solved using a Singular Value Decomposition. We now give a lower bound on the estimation accuracy. 
Let $\hat{R}$ be an unbiased estimator of $R$ in the sense of the intrinsic (right invariant) error $R \hat R^T$. After $n$ measurements, the information matrix is (using the independence of the noises)
$$
\xi^T J_n(R) \xi=-\sum_1^n \textbf{E}_R \left( \left. \frac{d^2}{dt^2} \right|_{t=0} \ln p \left( X_k\mid \exp(t\xi)R \right) \right).
$$
To derive $J$ we first note that:
$$
\norm{X_k-R^T\exp \left( - t \xi \right)d_k}^2 = \norm{X_k}^2-2d_k^T \exp \left( t \xi \right) RX_k+\norm{d_k}^2
$$
We have besides the following Taylor expansion 
\begin{align}\text{expm}\left( (t\xi)_\times \right)=I+t(\xi)_\times+\frac{1}{2}t^2 (\xi)_\times^2+O(t^3)\label{taylor:exp}\end{align}
so the second detivative at $t=0$ of $\norm{X_k-R^T\exp \left( - t \xi \right)d_k}^2$ is $- d_k^T  (\xi)_\times^2 R X_k$. As we have $\textbf{E}_R \left( R X_k \right)=d_k$ we obtain:
\begin{align*}
\textbf{E}_R \left( -d_k^T (\xi)_{\times}^2 RX  \right)   = -d_k^T (\xi)_{\times}^2 d_k  = -\xi^T (d_k)_{\times}^2 \xi.
\end{align*}
Let $H_k=(d_k)_\times$, we have proved the following relation:
$$
\textbf{E}_R \left( \left. \frac{d^2}{dt^2} \right|_{t=0} \ln p \left( X_k\mid \exp(t\xi)R \right) \right) = \frac{1}{\sigma_k^2}H_k^TH_k,
$$
giving immediatly:
\begin{prop}For Wahba's problem \cite{wahba1965}, the intrinsic covariance matrix defined by \eqref{P:def} where $\hat R$ is any unbiased estimator for a sample of $n$ independent measurements, satisfies inequality \eqref{result2} where
$$J(R)=\sum_1^n\frac{1}{\sigma_k^2}H_k^TH_k.$$
\end{prop} 
 
 Note that, the result does not depend on the underlying parameter $R$. This might be explained using theory of equivariant estimators on Lie groups that can be traced back to  \cite{shuster2006,peisakoff} (see also  \cite{berger1985statistical} for a more recent exposure).

\subsection{Proof of the result \eqref{result2}}
Let $\hat{R}$ be an unbiased estimator of $R$ in the sense of the intrinsic (right invariant) error $R \hat R^T$, that is,
\[
\int_X \log \left[ R \hat R^T (X) \right] p \left(X \mid R \right) dX =  0.
\]
So, if we let $\xi$ be any vector of the Lie algebra and $t\in\RR$, differentiating the latter equality we get
\[
\frac{d}{dt} \int_X \log \left[ \exp(t\xi)R\hat R ^T(X) \right] p \left( X \mid \exp(t\xi)R \right) dX =  0.
\]
Formally, this implies
\begin{equation}
\begin{aligned}
\label{eq1}
\int_X \Big( D\log & \left( R \hat{R}^T(X) \right) \left[ (\xi)_{\times} R \hat{R}^T(X)  \right] p(X|R) \\&+ \log (R \hat{R}^T(X) ) D_2p \left( X\mid R \right) [(\xi)_\times R] \Big) dX = 0\end{aligned},
\end{equation}
where $D$ denotes the differential and $D_2$ the partial differential with respect to the second argument. 
For any vector $u \in \frak{g}$ we have thus:
\begin{align*}
-\int_X  & \left\langle u,D\log \left( R \hat{R}^T(X) \right) \left[ (\xi)_{\times} R \hat{R}^T(X) \right] \right\rangle  p(X|R)dX \\
 & = \int_X \left\langle u, \log \left( R \hat{R}^T(X) \right) \right\rangle    D_2p(X\mid R) \left[ (\xi)_\times R \right] dX \\
 & \leqslant \sqrt{ \int_X \left\langle u, \log \left( R \hat{R}^T(X) \right) \right\rangle  ^2 p(X|R) dX} ~\cdot\\&\qquad \sqrt{\int_X \left( D_2 \ln p(X\mid R) [(\xi)_\times R] \right)^2 p(X|R) dX }
\end{align*}
where we used the Cauchy Schwarz inequality and the relationships
$$
D_2p=pD_2\log p,\quad \text{and then} \quad p=(\sqrt p)^2
$$
We introduce a basis of $\frak{g}$ and the matrix $\tilde{A}(X)=\log(R \hat{R}(X)^T)$. The latter inequality can be re-written as follows:
%\begin{equation}
\begin{align}
& \left[ u^T \int_X D\log \left( R \hat{R}^T(X) \right) \left[ (\xi)_{\times} R \hat{R}^T(X) \right] p(X|R)dX \right]^2  \nonumber \\ 
 & \leqslant \left( u^T \left[ \int_X \tilde{A}(X)^T \tilde{A}(X) p(X|R) dX \right] u \right) \xi^T J \left( R \right) \xi \label{Cauchy-Schwarz}
\end{align}
%\end{equation}
Now we compute a second-order expansion of the left-hand term  in the estimation error $Q:=R \hat{R}^T(X)$. To do so, we note 
 $
\log [ \exp(t\xi)Q  ]=\log [ \exp(t\xi)\exp(\log(Q))  ]$ is equal to (using the Baker-Campbell-Hausdorff formula and keeping only terms up to $t^2$)  $$\xi - \frac{1}{2} ad_{\log(Q)} t\xi+ \frac{1}{12} ad_{\log(Q)}^2  t\xi+ O(||\log(Q)||^3)] t\xi 
$$ Differentiating w.r.t to $t$ yields
\begin{align*}
D\log(Q) [(\xi)_{\times} Q ] = \Big[ I - \frac{1}{2} ad_{\log(Q)} +& \frac{1}{12} ad_{\log(Q)}^2  \\  &+ O \left(||\log(Q)||^3 \right) \Big] \xi 
\end{align*}
Neglecting the third-order terms we get:
\begin{align*}
&D\log(Q) [(\xi)_{\times} Q ] = \left[ I - \frac{1}{2} (\log(Q))_\times + \frac{1}{12} (\log(Q))_\times^2  \right] \xi =
\\& \left[ I - \frac{1}{2} (\log(Q))_\times + \frac{1}{12} \left( \log(Q) \log(Q)^T-||\log(Q)||^2 I_3 \right) \right] \xi
\end{align*}
We introduce the latter second-order expansion in the error in the equation \eqref{Cauchy-Schwarz} (note this is only formal as the second-order expansion will appear in an integral):
\begin{align*}
\bigg[u^T \int_X  &  \Big[ I - \frac{1}{2} (\log (R \hat{R}(X)^T))_\times + \frac{1}{12} (\log (R \hat{R}(X)^T \log (R \hat{R}(X)^T \\
- &\quad ||\log (R \hat{R}(X)^T)||^2 I_3) \Big] p(X|R) dX \xi \bigg]^2 \\
 & \leqslant \left(u^T \left[\int_X \tilde{A}(X)^T \tilde{A}(X) p(X|R) dX \right]u \right) \left( \xi^T J \left( R \right) \xi \right)
\end{align*}
This writes using the fact that $\hat R$ is unbiased (which makes the term in front of the factor $\frac{1}{2}$ cancel as $\textbf{E}_R \left( \log(Q) \right)=0$)
\begin{align*}
[u^T \xi + \frac{1}{12 } & u^T \int_X  [\tilde{A}(X) \tilde{A}(X)^T - Tr(\tilde{A}(X) \tilde{A}(X)^T) I_3] p(X|R) dX \xi ]^2 \\
 & \leqslant \left(u^T \left[\int_X \tilde{A}(X)^T \tilde{A}(X) p(X|R) dX \right] u \right) \left( \xi^T J \left( R \right) \xi \right)
\end{align*}
Letting $P=\int_X \tilde{A}(X) \tilde{A}(X)^T p(X|R) dX$ we get:
\[
\left[ u^T \left( \left[1 - \frac{1}{12} Tr(P) \right] I_3 + \frac{1}{12} P \right) \xi \right]^2 \leqslant (u^T P u) (\xi^T J \xi)
\]
The change of variables $\xi' = \left( \left[1 - \frac{1}{12} Tr(P) \right] I_3 + \frac{1}{12} P \right) \xi $ yields:
\begin{align*}
(u^T \xi')^2 \leqslant (u^T P u)~\cdot & \Bigg( {\xi'}^T \left( \left[1 - \frac{1}{12} Tr(P) \right] I_3 + \frac{1}{12} P \right)^{-1} \\& J \left( \left[1 - \frac{1}{12} Tr(P)  \right] I_3 + \frac{1}{12} P \right)^{-1} \xi' \Bigg)
\end{align*}
This inequality is true for any $u,\xi'$, which implies matricially the following desired result:
\begin{equation}
\label{result}
 P \succeq 
 \left( \left[1 - \frac{Tr(P)}{12}  \right] I_3 + \frac{1}{12} P \right) J^{-1}  \left( \left[1 - \frac{Tr(P)}{12} \right] I_3 + \frac{1}{12} P \right) 
\end{equation}
\subsection{Links with the Cram{\'e}r-Rao bound on manifolds \cite{smith-2005}}
The main result of \cite{smith-2005} stipulates that:
\begin{equation}
\label{smith}
P + \frac{1}{3} [R_m(P)J^{-1} + J^{-1} R_m(P) ] - \frac{1}{9} R_m(P) J^{-1} R_m(P)] \succeq J^{-1}
\end{equation}
where $R_m(P)$ is defined through Riemann's curvature tensor, and making use of the bi-invariance of the right-invariant metric  of $SO(3)$ (see e.g. \cite{arnold-fourier} for Riemann's sectional curvature formulas on Lie groups)  it boils down to $\langle R_m(P) u, u\rangle = \mathbb{E}(\frac{1}{4} || u \times \log (R^T \hat{R}(X) ) ||^2) = -\frac{1}{4} u^T \mathbb{E}((\log[R^T \hat{R}(X)])_\times^2) u = -\frac{1}{4} u^T[P - Tr(P)]u$. Replacing in \eqref{smith} we get:
\begin{align*}
P - &\frac{1}{12} \left([P-Tr(P)I_3] J^{-1} + J^{-1} [P-Tr(P)I_3] \right) \\&\qquad- \frac{1}{144}[P-Tr(P)I_3] J^{-1} [P-Tr(P)I_3] \succeq J^{-1}
\end{align*}
The same formula is obtained developing \eqref{result}, which proves both results coincide. However, the reader can check that the calculation on a general Riemannian manifold is  more tedious (and local). The Lie group case is more straithtforward and the results are obtained without using the machinery of Riemannian second-order geometry. 

\section{Application to dynamic attitude filtering}
\label{taylor:sec}

We now consider the following system
\begin{equation}
\begin{aligned}
\dotex R_t&=R_t(\omega(t))_\times, \quad
Y_n= \begin{pmatrix} Y_n^0 \\ Y_n^1 \end{pmatrix} = \begin{pmatrix} R_{t_n}^Td^0_n+V^0_n \\ R_{t_n}^Td^1_n+V^1_n \end{pmatrix} \label{deded}
\end{aligned}\end{equation}
that is, the motion in space of a solid fixed at a point, having deterministic known angular velocity  $\omega(t)$, and noisy measurements at discrete times $t_1\leq t_2\leq \cdots$. This fits into the general filtering setting of  \cite{taylor1979}. We assume that $V_n^0$ and $V_n^1$ are Gaussian noises with covariance matrices $\left( \sigma_n^0 \right)^2Id$ and $\left( \sigma_n^1 \right)^2Id$. Our goal is to derive an intrinsic lower Cram{\'e}r-Rao bound on the estimation error. We will see it follows from the results of Section \ref{smith:sec} indeed, thanks to the facts that 1- $R_t$ is a deterministic quantity and 2- due to the invariance of the system, the flow can be explicitly computed. 

Such problems arise for attitude estimation in the degenerate case where the gyroscope is infinitely better than the vector sensors. Sensors measuring in the body frame vectors from the fixed frame include magnetometers, that measure the earth magnetic field in the body frame, and accelerometers, that measure the earth gravity vector field in the body frame, under static flight assumptions. For each of these sensors, the isotropy assumption of the noise is reasonable technologically, as the measurements are performed using in each case three orthogonal one-axis sensors (accelerometers or magnetometers). The Gaussianity is more questionable but it is a convenient and widespread assumption about the noise. 

\subsection{Intrinsic Fisher information matrix  computation}
The conditional intrinsic information matrix at time $n$ is (using \cite{taylor1979} and the results above)
$$\xi^TJ_n\xi=- \textbf{E}_R \left( \frac{d^2}{dt^2}\ln p \left(Y_1,\cdots, Y_n\mid \exp(t\xi)R_{t_n} \right) \right).$$
 Now, using the invariance of the dynamics, we see there exists a rotation $A(k,n)$ depending only on $\omega(s),~t_k\leq s\leq t_n$ such that $R_{t_n}=R_{t_k}A(k,n)$. Indeed $A(k,n)$ is the solution at time $t_n$ to the differential equation on $SO(3)$ defined by $\dotex A_s=A_s\omega(s),~A_{t_k}=R_{t_k}$ (see eg.  \cite{barrau2013intrinsicp}). As a result, all the measurements are independent given $R_{t_n}$ and we can write:
\begin{align*}
-\ln p & \left(Y_1,\cdots, Y_n\mid \exp(t\xi)R_{t_n} \right) \\
 = & - \sum_{k=1}^n \ln p \left(Y_k \mid \exp(t\xi) R_{t_n} A(k,n)^T \right) \\
 = & Cst + \sum_{k=1}^n\frac{1}{2\sigma_k^2}\norm{Y^0_k-A(k,n)R_{t_n}^T\exp(-t\xi)d^0_k}^2 \\
 & + \sum_{k=1}^n\frac{1}{2\sigma_k^2}\norm{Y^1_k-A(k,n)R_{t_n}^T\exp(-t\xi)d^1_k}^2
\end{align*}
When deriving an intrinsic Cram{\'e}r-Rao bound for Wahba's problem in Section \ref{smith:sec}, we have already proved that letting $H=(d)_\times$ we have for any $Q\in SO(3)$ such that $ \textbf{E}(QY)=d$ that \begin{equation}
\begin{aligned}
\frac{d^2}{dt^2} \textbf{E} \left( \norm{Y-Q^T\exp(-t\xi)d}^2 \right)&=-2\frac{d^2}{dt^2} \textbf{E} (d^TQ\exp(-t\xi)Y)\\&=2\xi^THH^T\xi\label{ee:eq}
\end{aligned}\end{equation}
and the result does not depend on $Q$. Differentiating $-\ln p \left(Y_1,\cdots, Y_n\mid \exp(t\xi)R_{t_n} \right)$ twice w.r.t. $t$ and using \eqref{ee:eq} with $Q=R_{t_n}A(k,n)^T$ (which is valid as $\textbf{E}_{R_{t_k}} (QY_k^i)= \textbf{E}_{R_{t_k}} \left( R_{t_n}A(k,n)^TY_k^i \right)= \textbf{E}_{R_{t_k}} (R_{t_k}Y_k^i)=d_k^i$ for $i=1,2$), we get finally:
\begin{align*}
J_n&=-\frac{d^2}{dt^2} \ln p(Y_1,\cdots, Y_n\mid R_{t_1},\cdots ,\exp(t\xi)R_{t_n})\\& = \sum_{k=1}^n \left( \frac{1}{\left( \sigma_k^0 \right)^2} {H_k^0}^T H_k^0 + \frac{1}{\left( \sigma_k^1 \right)^2} {H_k^1}^T H_k^1 \right)
\end{align*}

\subsection{Intrinsic Cram{\'e}r-Rao bound }A mere application of Theorem 1 implies
\begin{prop}For the considered system \eqref{deded}, at time $t_n$, the accuracy $P=E(\log(R_{t_n}\hat R_{t_n}^T)\log(R_{t_n}\hat R_{t_n}^T)^T)$ of any unbiased attitude estimator is lower bounded according to formula \eqref{result2}, with  $J(R)= \sum_{k=1}^n \left( \frac{1}{\left( \sigma_k^0 \right)^2} {H_k^0}^T H_k^0 + \frac{1}{\left( \sigma_k^1 \right)^2} {H_k^1}^T H_k^1 \right)$. \end{prop}

\section{Links with Invariant Kalman filtering}\label{kalman:sec}

For the filtering problem of Section \ref{taylor:sec}, on can derive an Invariant Extended Kalman Filter (IEKF) \cite{bonnabel2007left}. The IEKF is a novel methodology for devising EKFs on Lie groups, where the EKF is bound to respect the invariances of the problem, and where an intrinsic estimation error is linearized at each step. Moreover, the exponential map allows to map the Kalman correction term to the state space. The IEKF for the problem above is derived on $SO(3)$ in the recent paper \cite{barrau2013intrinsicp}, and we briefly  recall the principle here. The IEKF equations write
\begin{align}
\dotex\hat R_t&=\hat R_t(\omega(t))_\times,\quad t_{n-1}<t<t_{n}\qquad\text{(Propagation)}\\ \hat R_{t_n}^+&=\exp \left(K_n  \begin{pmatrix} \hat R_{t_n} Y^0_n-d^0_n \\ \hat R_{t_n} Y^1_n-d^1_n \end{pmatrix} \right)\hat R_{t_n},\quad t=t_n \qquad \text{(Update)}
\end{align}
where $K_n\in\RR^{3\times 3}$ is the gain matrix to be tuned as follows.  Letting $\xi_t=\log \left( \hat R_t R_t^T \right)$ be the right invariant estimation error projected in the Lie algebra, the error system has the following remarkable autonomous form
\begin{align*}
\dotex\xi_t&=0,\quad t_{n-1}<t<t_{n}\qquad\text{(Propagation)} \\
\exp(\xi_{t_n})^+&=\exp \left( K_n \begin{pmatrix} \exp(\xi_{t_n}) d_n^0-d_n^0+V_n^0 \\
\exp(\xi_{t_n}) d_n^1-d_n^1+V_n^1 \end{pmatrix} \right) \exp(\xi_{t_n}) \\
 & \tag{ Update}
\end{align*}
During the propagation step, the covariance of the linearized estimation error $P_t=\textbf{E} \left( \xi_t\xi_t^T \right)$ remains fixed, that is,
$$P_{t_{n+1}}=P_{t_n}^+$$
as the linearized dynamics (for the well-chosen estimation error) yields a static system and it was assumed there is no process noise. As concerns the update step, using formula \eqref{taylor:exp}, a \emph{first order approximation} to the error update equation above reads \cite{barrau2013intrinsicp}
\begin{align*}
\xi_{t_n}^+
& = \xi_{t_n}+  K_n \begin{pmatrix} \xi_{t_n}\times d^0_n + V_n^0\\\xi_{t_n}\times d^1_n + V_n^1 \end{pmatrix} \\
 & = \left( I_3-K_n \begin{pmatrix} H_n^0 \\ H_n^1 \end{pmatrix} \right) \xi_{t_n} +K_n \begin{pmatrix} V_n^0 \\ V_n^1 \end{pmatrix}.
\end{align*}
The gain that minimizes the increase in the covariance matrix of the linearized error at the update step is the Kalman gain
$$
K_n=P_{t_n}\begin{pmatrix} H_n^0 \\ H_n^1 \end{pmatrix}^T \left( \begin{pmatrix} H_n^0 \\ H_n^1 \end{pmatrix} P_{t_n} \begin{pmatrix} H_n^0 \\ H_n^1 \end{pmatrix}^T + N  \right)^{-1},
$$
with $N=\begin{pmatrix} \left( \sigma_n^0 \right)^2 I_3 & 0_3 \\ 0_3 &  \left( \sigma_n^1 \right)^2 I_3 \end{pmatrix}$,
leading to the covariance update:
$$
(P_{t_n}^+)^{-1}=P_{t_n}^{-1} +   \frac{\left( H_n^0 \right)^T \left( H_n^0 \right)}{ \left( \sigma_n^0 \right)^2}   + \frac{\left( H_n^1 \right)^T \left( H_n^1 \right)}{ \left( \sigma_n^1 \right)^2}  
$$
As there is no a priori information about the value of $R_{t_1}$, a usual way to initialize the filter is maximum likelihood:
\begin{equation}
\label{eq::max_like}
\hat R_{t_1} = \min_{\hat R} \left \{ \norm{\hat R Y_1^0 - d_1^0}^2/\left( \sigma_1^0 \right)^2 + \norm{ \hat R Y_1^1 - d_1^1 }^2/\left( \sigma_1^1 \right)^2 \right \}
\end{equation}
As concerns the covariance matrix $P_1= Cov \left( \xi_{t_1} \right)$, where $\exp \left( \xi_{t_1} \right)= \hat R_{t_1} R_{t_1}^T$, a first-order expansion of  \eqref{eq::max_like} reads $\xi_{t_1} = \min \norm{\xi_{t_1} \times d_1^0 + V_n^0}^2 +  \norm{\xi_{t_1} \times d_1^1 + V_n^1}^2$, i.e. $\xi_{t_1} = \left( \frac{\left( H_1^0 \right)^T \left( H_1^0 \right)}{ \left( \sigma_1^0 \right)^2}   + \frac{\left( H_1^1 \right)^T \left( H_1^1 \right)}{ \left( \sigma_1^1 \right)^2}    \right)^{-1}    \left( H_1^0 V_1^0 + H_1^1 V_1^1 \right)
$ which gives: $ \left( P_{t_1} \right)^{-1} = Cov \left( \xi_{t_1} \right)^{-1} =  \frac{\left( H_1^0 \right)^T \left( H_1^0 \right)}{ \left( \sigma_1^0 \right)^2}  +  \frac{\left( H_1^1 \right)^T  \left( H_1^1 \right)}{ \left( \sigma_1^1 \right)^2} $. Gathering the previous results we obtain:

\begin{prop}The covariance matrix of the error returned by the IEKF writes
$$
P_n=
\left( \sum_{k=1}^n\frac{{H_k^0}^T H_k^0}{ \left( \sigma_k^0\right)^2} + \sum_{k=1}^n\frac{{H_k^1}^T H_k^1}{ \left( \sigma_k^1\right)^2} \right)^{-1}
$$
and thus the IEKF returns the Cram{\'e}r-Rao bound for the associated filtering problem, neglecting the curvature terms.\end{prop} 

Note that,  it is logical that the  curvature terms be ignored by the IEKF as it is based on a first order approximation of the estimation error.
The result is in sharp contrast with the general theory \cite{taylor1979} that stipulates that the Cram{\'e}r-Rao bound is the EKF covariance indeed, but, linearized around the true trajectory, that is unknown to the user. Those bounds are referred to as ``posterior Cram{\'e}r-Rao bounds'' in the filtering and hidden Markov models (HMM) literature (see e.g. \cite{tichavsky1998posterior}). For invariant systems on SO(3), in the case of deterministic dynamics, we have proved the bound can be computed in real time.  This appears as another remarkable feature of dynamical systems defined on Lie groups.

\section{Conclusion}
In this paper we have derived an intrinsic Cram{\'e}r-Rao lower bound (ICRLB) on $SO(3)$ in a straightforward way. We have applied it to derive an ICRLB for Wahba's problem when the noise is isotropic and Gaussian. Then, we have also derived an ICRLB for the problem of filtering on $SO(3)$ a system with deterministic evolution and noisy isotropic Gaussian measurements.  We have also proved the \emph{intrinsic} CRLB is the covariance matrix returned by the \emph{invariant} EKF on $SO(3)$. This is a remarkable result, as generally the CRLB can not be computed online as it presupposes to know the true trajectory of the system, which is precisely what one seeks to estimate. It is thus usually reserved for offline simulations to test filters' efficiency, and we generally speak of ``posterior Cram{\'e}r-Rao bounds''   \cite{tichavsky1998posterior}. 

In the future, we would like to investigate  in what ways the intrinsic gradient methods (see \cite{bonnabel2013stochastic,amari-98}) might asymptotically reach the Intrinsic Cram{\'e}r-Rao bound. Besides, we hope it is possible to derive an ICRLB for the filtering problem considered in the present paper, but with a noisy evolution (that is, using noisy gyroscopes). In this case it will certainly not coincide with the covariance returned by the IEKF, but it might be computable online taking advantage of the invariances of the system. Such results could be applied to a wide range of aeronautics estimation problems, like e.g., \cite{barczyk2013invariant}, and may be useful to test efficiency of some other instrinsic filtering methods, such as \cite{bourmaud2013discrete}.

%\newtheorem{theorem}{Proposition}

%\subsection{Considered continuous-time model}

\bibliographystyle{plain}

%\bibliographystyle{IEEEtran}
%\bibliography{IEEEabrv,rhn}

\end{document}